# EVALUATION OF FORMAL POSTERIOR DISTRIBUTIONS VIA MARKOV CHAIN ARGUMENTS

By Morris L. Eaton, James P. Hobert,[1]
Galin L. Jones and Wen-Lin Lai

*University of Minnesota, University of Florida, University of Minnesota and Providence University*

We consider evaluation of proper posterior distributions obtained from improper prior distributions. Our context is estimating a bounded function $\phi$ of a parameter when the loss is quadratic. If the posterior mean of $\phi$ is admissible for all bounded $\phi$, the posterior is *strongly admissible*. We give sufficient conditions for strong admissibility. These conditions involve the recurrence of a Markov chain associated with the estimation problem. We develop general sufficient conditions for recurrence of general state space Markov chains that are also of independent interest. Our main example concerns the $p$-dimensional multivariate normal distribution with mean vector $\theta$ when the prior distribution has the form $g(\|\theta\|^2)\,d\theta$ on the parameter space $\mathbb{R}^p$. Conditions on $g$ for strong admissibility of the posterior are provided.

**1. Introduction.** In many standard parametric settings, the use of improper prior distributions to produce inferential proposals is well established. A variety of formal rules have been proposed to justify particular improper priors. This is especially true in invariant problems. The survey article of Kass and Wasserman [15] provides an excellent overview of many issues that arise in selecting and using improper priors, including decision theoretic considerations. Indeed, decision theory provides an appealing framework for the evaluation of improper priors via the (proper) posterior distributions they produce. A particular decision theoretic approach was described in [7], but subsequent suggestions have had more appeal. One of these involves a notion of "strong admissibility" and is the basis of the approach below. This notion was introduced in [8] and was called $P$-admissibility in [14].

Received July 2007; revised July 2007.
[1]Supported in part by NSF Grant DMS-05-03648.

The idea behind strong admissibility is the following. Given a sample space $(\mathcal{X}, \mathcal{B})$, suppose a random object $X \in \mathcal{X}$ is to be observed and assume that the distribution of $X$ is an element of the parametric model $\{P(\cdot|\theta)|\theta \in \Theta\}$. Consider a $\sigma$-finite improper prior distribution $\nu$ and define the marginal measure $M$ on $(\mathcal{X}, \mathcal{B})$ by

$$M(B) = \int_\Theta P(B|\theta)\nu(d\theta), \quad B \in \mathcal{B}.$$

Technical issues such as measurability, existence of integrals, etc. are treated carefully in the next section. When the measure $M$ is $\sigma$-finite, a proper posterior distribution $Q(d\theta|x)$ exists and satisfies

$$P(dx|\theta)\nu(d\theta) = Q(d\theta|x)M(dx),$$

in the sense that the joint measures on each side of the equality sign are equal. Interpreting $Q(\cdot|x)$ as summarizing knowledge about $\theta$ after seeing $X = x$, the probability measure $Q(\cdot|x)$ on $\Theta$ can now be used to solve decision problems. Namely, one integrates a loss function with respect to $Q(\cdot|x)$ and then picks an action (depending on $x$) to minimize the expected posterior loss. This procedure is ordinarily called the formal Bayes method of solving a decision problem.

As an example, let $\phi(\theta)$ be a bounded real-valued function on $\Theta$ and consider a quadratic loss function

$$(1) \qquad L(a, \theta) = (a - \phi(\theta))^2, \qquad a \in \mathbb{R}.$$

The posterior expected loss is minimized by the posterior mean

$$(2) \qquad \hat{\phi}(x) = \int_\Theta \phi(\theta) Q(d\theta|x),$$

and $\hat{\phi}$ is a formal Bayes estimator of $\phi$. Given any estimator $t(x)$ of $\phi(\theta)$, the *risk function* of $t$ is

$$(3) \qquad r(t, \theta) = \int_\mathcal{X} (t(x) - \phi(\theta))^2 P(dx|\theta),$$

and admissibility of estimators is assessed in terms of the risk function.

The formal posterior $Q(d\theta|x)$ is *strongly admissible* if, for every bounded measurable function $\phi$, the formal Bayes estimator $\hat{\phi}$ is admissible (this is defined more carefully in the next section). Loosely speaking, strong admissibility is regarded as an endorsement of $Q(\cdot|x)$ and, hence, of the improper prior $\nu$, for use in making inferences about $\theta$ after seeing $X = x$.

A sufficient condition for strong admissibility is the recurrence of the Markov chain with state space $\Theta$ and transition function

$$(4) \qquad R(C|\theta) = \int_\mathcal{X} Q(C|x) P(dx|\theta).$$

In words, (4) is the expectation of the formal posterior when $X$ is sampled from $P(\cdot|\theta)$. There is more than one notion of recurrence that is relevant for strong admissibility. A discussion of this issue is given in the next section. Here is an example that illustrates the issues that motivated our research.

EXAMPLE 1.1. Suppose $X$ has a $p$-dimensional multivariate normal distribution with mean vector $\theta \in \mathbb{R}^p$ and covariance matrix the $p \times p$ identity, $I_p$. This we write as $X \sim N_p(\theta, I_p)$. Consider Lebesgue measure as the improper prior distribution: $\nu(d\theta) = d\theta$ on $\Theta = \mathbb{R}^p$. Standard calculations show that the formal posterior $Q(d\theta|x)$ is a $N_p(x, I_p)$ distribution on $\mathbb{R}^p$. Further, the transition function $R(\cdot|\theta)$ is a $N_p(\theta, 2I_p)$ distribution on $\mathbb{R}^p$. Therefore, the one step transition of the Markov chain in this example can be described as follows: Given the chain is at $\theta \in \mathbb{R}^p$, the next state of the chain is $\theta + V$, where $V$ is $N_p(0, 2I_p)$. The chain is thus a random walk on $\mathbb{R}^p$. For $p = 1$ or $p = 2$, the chain is recurrent (see [5] for $p=1$ and [20], Chapter 3, for $p=2$). But for $p \geq 3$, the chain is transient ([12], page 579). Therefore, for $p = 1, 2$, the formal posterior $Q(\cdot|x)$ is strongly admissible, but for $p \geq 3$, the recurrence argument fails. Of course, one suspects that, for $p \geq 3$, $Q(\cdot|x)$ is in fact not strongly admissible because of the existence of the James–Stein estimator for $\theta$, but a rigorous proof of a "not strongly admissible" assertion for $Q(\cdot|x)$ is not known to us.

Now, focus on the case of $p \geq 3$. Since the chain induced by the normal model and the improper prior "$d\theta$" is transient, it seems reasonable to look at a somewhat broader class of priors of the form

$$\nu(d\theta) = g_0(\|\theta\|^2) \, d\theta, \tag{5}$$

where $\|\cdot\|$ denotes the Euclidean norm on $\mathbb{R}^p$. It is natural to seek conditions on the function $g_0$ so that the induced Markov chain is recurrent. It is this and related issues that give rise to much of the new material herein. The primary difficulties in addressing this issue are as follows:

(i) The transition function cannot be computed explicitly, so the corresponding Markov chain is difficult to study.

(ii) Even for special and rather simple $g_0$'s, techniques for establishing the recurrence of the induced chain are not available.

This completes our initial discussion of this example.

Here is an outline of the material in this paper. The next section contains notation and assumptions, a careful definition of strong admissibility, and a discussion of the Markov chain introduced above. In particular, two different notions of recurrence are described (neither of which involves the notion of irreducibility), and their relation to strong admissibility is emphasized. Results in [8, 9, 10] underlie much of this material.

Section 3 contains the first of two main results that are of use in establishing strong admissibility. This result was conjectured in [8], pages 1166–1167. The idea, reminiscent of hierarchical prior proposals, is to give conditions under which the improper prior distribution $\nu(d\theta)$ on $\Theta$ has a "disintegration" into $\pi(d\theta|\beta)s(d\beta)$. Here $\pi(\cdot|\beta)$ is a probability measure for each $\beta$ and $s$ is a $\sigma$-finite measure on a set $A$ (typically, much simpler than $\Theta$). Sufficient conditions for strong admissibility are then based on an induced Markov chain whose state space is $A$ rather than $\Theta$. The arguments in Section 3 rely on results involving the Dirichlet forms associated with the Markov chains under consideration (see [10]).

Our second main result, presented in Section 4, concerns sufficient conditions for the recurrence of a Markov chain with a general state space (a Polish space) $Y$. The results in Section 4 are an extension and refinement of results in [16], which are in turn an extension of the early work in [17]. When $Y = [0, \infty)$, sufficient conditions for the recurrence of a set $[0, m)$ for $m$ sufficiently large are expressed in terms of the first three moments of the increments of the chain.

In brief, recurrence of the Markov chain with transition function $R$ implies strong admissibility. This paper provides two results that simplify the task of establishing recurrence. Application of the results to Example 1.1 illustrates a nontrivial case that justifies the new methodology in Sections 3 and 4. In particular, these methods are applied in Section 5 to the problem of Example 1.1 when the prior has the form given in (5). Conditions on the function $g_0$ (which are necessarily dimension dependent) that imply strong admissibility are discussed in detail. Our treatment of Example 1.1 is a rigorous development of work begun in [16].

**2. Background.** The purpose of this section is to provide assumptions and rigorous statements that underlie the connection between strong admissibility and recurrence. The discussion here is an abbreviated version of some material in [10] and we assume the reader is somewhat familiar with that material. It is useful to keep Example 1.1 in mind.

2.1. *Model, prior and posterior.* The sample space $(\mathcal{X}, \mathcal{B})$ is assumed to consist of a Polish space $\mathcal{X}$ (a complete separable metric space) coupled with the Borel $\sigma$-algebra $\mathcal{B}$. The parameter space $(\Theta, \mathcal{C})$ is also assumed to be Polish with $\mathcal{C}$ the Borel $\sigma$-algebra. A statistical model $\{P(\cdot|\theta)|\theta \in \Theta\}$ for an observable quantity $X \in \mathcal{X}$ specifies the modeling assumption. We assume the $P(\cdot|\cdot)$ are Markov transition functions—that is, $P(\cdot|\theta)$ is a probability measure on $\mathcal{B}$ for each $\theta$, and for each $B \in \mathcal{B}$, $P(B|\cdot)$ is $\mathcal{C}$-measurable.

Now let $\nu(d\theta)$ be a $\sigma$-finite measure on $\mathcal{C}$, and for each $B \in \mathcal{B}$, consider the marginal measure given by

$$(6) \qquad M(B) = \int_\Theta P(B|\theta)\nu(d\theta).$$

Throughout, the marginal measure $M$ is assumed to be $\sigma$-finite. Under this assumption, there is a Markov transition function $Q(\cdot|x)$ that satisfies

(7) $$P(dx|\theta)\nu(d\theta) = Q(d\theta|x)M(dx).$$

This equation means that the two joint measures on $\mathcal{X} \times \Theta$ given by (7) agree. For a discussion of the existence and uniqueness (a.e. $M$) of $Q(\cdot|x)$; see [7]. Of course, $Q(\cdot|x)$ is called the *formal posterior distribution* of $\theta$ given $X = x$.

2.2. *Strong admissibility and recurrence.* As in Section 1, let $\phi$ be a bounded, measurable, real-valued function defined on $\Theta$ and consider the problem of estimating $\phi(\theta)$ when the loss is (1). Then $\hat{\phi}(x)$ given in (2) is the formal Bayes estimator obtained from the formal posterior $Q(\cdot|x)$. Using the risk function defined in (3), here is an appropriate notion of admissibility, due to C. Stein, for our setting.

DEFINITION 2.1. The estimator $t_0(x)$ for $\phi(\theta)$ is *almost-$\nu$-admissible* (a-$\nu$-a) if for each estimator $t$ that satisfies $r(t,\theta) \leq r(t_0,\theta)$ for all $\theta$, the set $\{\theta | r(t,\theta) < r(t_0,\theta)\}$ has $\nu$-measure zero.

DEFINITION 2.2. The improper prior $\nu$, or equivalently, the formal posterior $Q(\cdot|x)$, is *strongly admissible* if for each bounded measurable $\phi$, the estimator $\hat{\phi}$ is a-$\nu$-a.

The boundedness assumption on $\phi$ greatly simplifies the technical issues surrounding our discussion. For some parallel results regarding the estimation of unbounded functions, see [9].

In this section two notions of recurrence are useful. To describe these, recall the transition function given in (4). This transition function determines a Markov chain $W = (W_0, W_1, \ldots)$ with each $W_i \in \Theta$, $i = 0, 1, 2, \ldots$. The path space of $W$ is $\Theta^\infty$ and given that $W_0 = w_0$, the probability distribution of $W$ is denoted by $P_{w_0}$. In some situations we will call $W$ the $P$-$\nu$ chain to emphasize its dependence on the model and the prior.

A measurable subset $C \subseteq \Theta$ is $\nu$-*proper* if $0 < \nu(C) < \infty$. Let $C$ be a $\nu$-proper set and consider the stopping time $\tau_C$ defined on $\Theta^\infty$ by

$$\tau_C = \begin{cases} \infty, & \text{if } W_n \notin C \text{ for all } n \geq 1, \\ \text{smallest } n \geq 1 \text{ such that } W_n \in C, & \text{otherwise.} \end{cases}$$

Also, let $E_C = \{w | \tau_C(w) < \infty\} \subseteq \Theta^\infty$.

DEFINITION 2.3. The $\nu$-proper set $C$ is *locally-$\nu$-recurrent* (l-$\nu$-r) if the set $\{w_0 \in C | P_{w_0}(E_C) < 1\}$ has $\nu$-measure zero. The $\nu$-proper set $C$ is $\nu$-*recurrent* if the set $\{w_0 | P_{w_0}(E_C) < 1\}$ has $\nu$-measure zero.

DEFINITION 2.4. The Markov chain $W$ is *locally-$\nu$-recurrent* (l-$\nu$-r) if every $\nu$-proper set is l-$\nu$-r.

The following basic result was established in [8].

THEOREM 2.1. *If the Markov chain $W$ is l-$\nu$-r, then $Q(\cdot|x)$ is strongly admissible.*

The above theorem is often difficult to apply since every $\nu$-proper $C$ must be shown to be l-$\nu$-r. The following result ([10], Theorem 3.1), eases the burden somewhat.

THEOREM 2.2. *The following are equivalent:* (i) *the chain $W$ is l-$\nu$-r and* (ii) *there is an increasing sequence of $\nu$-proper sets $C_i$, $i = 1, 2, \ldots$, such that $C_i \to \Theta$ and each $C_i$ is l-$\nu$-r.*

The next result facilitates the application of the recurrence ideas to the admissibility issue. It is sometimes called the "one-set criterion." See [10] for a proof.

THEOREM 2.3. *Let $C^*$ be a $\nu$-proper set and suppose $C^*$ is $\nu$-recurrent. Then every $\nu$-proper set $C$ is l-$\nu$-r.*

Section 4 is devoted to the problem of finding a single recurrent set so that the above theorem can be applied to nontrivial examples.

The Markov chain one-set technique of Section 4 and the dimension reduction method of Section 3 are the basic contributions of this paper. The dimension reduction method relies on a rather deep connection between l-$\nu$-r and the behavior of the Dirichlet form associated with the Markov chain $W$ generated by the transition function $R(d\theta|\eta)$ and the measure $\nu$.

In what follows, $L_2(\nu)$ denotes the linear space of $\nu$-square integrable functions. Let $h \in L_2(\nu)$, then the quadratic form

$$(8) \qquad \Delta(h) = \tfrac{1}{2} \int \int (h(\theta) - h(\eta))^2 R(d\theta|\eta)\nu(d\eta)$$

is a Dirichlet form (see [6] and [13] for background material). The relevance of $\Delta$ for statistical decision problems with quadratic loss is discussed in [10]. Given a $\nu$-proper set $C$, let $I_C$ denote the indicator function of $C$ and let

$$(9) \qquad V(C) = \{h \mid h \in L_2(\nu), h \text{ bounded}, h \geq I_C\}.$$

Here is a basic theorem established in [8] that relates the Dirichlet form $\Delta$ to questions regarding admissibility.

THEOREM 2.4. *Fix a $\nu$-proper set $C$. The following are equivalent:* (i) *$C$ is l-$\nu$-r and* (ii) $\inf_{h \in V(C)} \Delta(h) = 0$.

Theorems 2.2, 2.3 and 2.4 are used in the next section.

**3. A reduction argument.** Again consider a model $\{P(dx|\theta) \mid \theta \in \Theta\}$, a $\sigma$-finite prior $\nu(d\theta)$ with $M(dx)$ assumed to be $\sigma$-finite. Hence, there is a Markov kernel $Q(d\theta|x)$, a formal posterior for $\theta$ given $x$, so that (7) holds. The transition function $R(d\theta|\eta)$ defines the $P$-$\nu$ Markov chain $W$ and the related symmetric measure

$$R(d\theta|\eta)\nu(d\eta) = \int_{\mathcal{X}} Q(d\theta|x)Q(d\eta|x)M(dx).$$

This in turn defines the Dirichlet form $\Delta$ given in (8). Next, consider a measurable mapping $t$ from $(\Theta, \mathcal{C})$ to $([0, \infty), \mathcal{B}_1)$. The map $t$ induces a measure $s$ on the Borel subsets $\mathcal{B}_1$ of $[0, \infty)$ given by $s(A) = \nu(t^{-1}(A))$.

ASSUMPTION 3.1. There is a sequence of disjoint Borel sets $A_i \subseteq [0, \infty)$ for $i = 1, 2, 3, \ldots$ such that

$$\bigcup_{i=1}^{\infty} A_i = [0, \infty) \quad \text{and} \quad B_i = t^{-1}(A_i) \subseteq \Theta \text{ is } \nu\text{-proper.}$$

THEOREM 3.1. *Suppose Assumption 3.1 holds. Then there is a Markov transition function $\pi(d\theta|a)$ on $\mathcal{C} \times [0, \infty)$ such that*

(10) $$\nu(d\theta) = \pi(d\theta|a)s(da).$$

*Equation* (10) *means that for all measurable nonnegative functions $f_1(\theta)$ and $f_2(t(\theta))$ defined on $\Theta$ and $[0, \infty)$ respectively, we have*

(11) $$\int_\Theta f_2(t(\theta))f_1(\theta)\nu(d\theta) = \int_0^\infty f_2(a)\left(\int_\Theta f_1(\theta)\pi(d\theta|a)\right)s(da).$$

If $\nu$ is a finite measure, then the assertion of Theorem 3.1 is nothing more than the existence of a conditional distribution [on $\Theta$ given $t(\theta) = a$] $\pi(d\theta|a)$, where $s$ is the marginal distribution of $t(\theta)$. When $\nu$ is allowed to be $\sigma$-finite, Assumption 3.1 allows one to argue separately on the spaces $B_i \subseteq \Theta$ and $A_i \subseteq [0, \infty)$ using standard techniques; see [18], pages 43 to 52. The details are omitted.

EXAMPLE 3.1. Consider $\nu(d\theta) = g_0(\|\theta\|^2) d\theta$ on $\Theta = \mathbb{R}^p$ with $p > 1$. Assume the measurable function $g_0$ is non-negative and $\nu$ is $\sigma$-finite. Consider the mapping $t$ on $\Theta$ defined by $t(\theta) = \|\theta\|^2 \in [0, \infty)$. Let $\pi(d\theta|\beta)$ denote the uniform probability distribution on $\{\theta|\|\theta\|^2 = \beta\}$ with the obvious extension of $\pi(\cdot|\beta)$ to $\Theta$. A routine argument shows that $\nu(d\theta) = \pi(d\theta|\beta)s(d\beta)$, where

$$s(d\beta) = \frac{[\Gamma(1/2)]^p}{\Gamma(p/2)} g_0(\beta)\beta^{p/2-1}\, d\beta \qquad \text{on } [0,\infty). \tag{12}$$

REMARK 3.1. Without Assumption 3.1, Theorem 3.1 is typically false. To see part of the difficulty, let $\Theta = [0,\infty) \times [0,\infty)$, let $\nu$ be Lebesgue measure on $\Theta$ and let $t(\theta_1, \theta_2) = \theta_1 \in [0,\infty)$. It is not too hard to show that a $\pi(d\theta|a)$ satisfying (11) does not exist.

For the remainder of this section, Assumption 3.1 holds, so the conclusion of Theorem 3.1 holds. Now, introduce the new model

$$\widetilde{P}(dx|a) = \int_\Theta P(dx|\theta)\pi(d\theta|a), \qquad a \in [0,\infty).$$

The "prior" $s(da)$ is $\sigma$-finite (this follows from Assumption 3.1) and the marginal on $\mathcal{X}$ is the $\sigma$-finite measure $M(dx)$. That $M(dx)$ is still the marginal follows from

$$\int_0^\infty \widetilde{P}(dx|a)s(da) = \int_0^\infty \int_\Theta P(dx|\theta)\pi(d\theta|a)s(da)$$
$$= \int_\Theta P(dx|\theta)\nu(d\theta) = M(dx).$$

Thus, there is a formal posterior $\widetilde{Q}(da|x)$ that satisfies

$$\widetilde{P}(dx|a)s(da) = \widetilde{Q}(da|x)M(dx). \tag{13}$$

The model $\widetilde{P}(dx|a)$ together with $s$ determines a Markov chain, the $\widetilde{P}$-$s$ chain, with transition function

$$\widetilde{R}(da|b) = \int \widetilde{Q}(da|x)\widetilde{P}(dx|b)$$

and Dirichlet form

$$\widetilde{\Delta}(\widetilde{h}) = \tfrac{1}{2}\int_0^\infty \int_0^\infty (\widetilde{h}(a) - \widetilde{h}(b))^2 \widetilde{R}(da|b)s(db).$$

The rest of this section is devoted to proving that if the $\widetilde{P}$-$s$ chain is locally-$s$-recurrent (l-$s$-r), then the $P$-$\nu$ chain is l-$\nu$-r. This shows that the l-$s$-r of the $\widetilde{P}$-$s$ chain implies strong admissibility for the original problem (see Theorem 2.1). To describe our first result, let $\widetilde{Q}_0$ be the Markov kernel on $\mathcal{B}_1 \times \mathcal{X}$ defined by

$$\widetilde{Q}_0(A|x) = Q(t^{-1}(A)|x). \tag{14}$$

THEOREM 3.2. *The Markov kernel $\widetilde{Q}_0$ serves as a version of $\widetilde{Q}$ in* (13).

PROOF. Consider nonnegative measurable functions $v$ on $\mathcal{X}$ and $\psi$ on $[0, \infty)$. Then using the definition of $\widetilde{Q}_0$ yields

$$\int_0^\infty \int_{\mathcal{X}} v(x)\psi(a)\widetilde{P}(dx|a)s(da) = \int_0^\infty \int_{\mathcal{X}} \int_\Theta v(x)\psi(a)P(dx|\theta)\pi(d\theta|a)s(da)$$
$$= \int_\Theta \int_{\mathcal{X}} v(x)\psi(t(\theta))P(dx|\theta)\nu(d\theta)$$
$$= \int_{\mathcal{X}} \int_\Theta v(x)\psi(t(\theta))Q(d\theta|x)M(dx)$$
$$= \int_{\mathcal{X}} \int_0^\infty v(x)\psi(a)\widetilde{Q}_0(da|x)M(dx). \quad \square$$

THEOREM 3.3. *Consider $\widetilde{h} \in L_2(s)$ and define $h^*$ on $\Theta$ by $h^*(\theta) = \widetilde{h}(t(\theta))$. Then $h^* \in L_2(\nu)$ and $\widetilde{\Delta}(\widetilde{h}) = \Delta(h^*)$, where $\widetilde{\Delta}$ and $\Delta$ are the Dirichlet forms for the $\widetilde{P}$-$s$ chain and the $P$-$\nu$ chain respectively.*

PROOF. That $h^* \in L_2(\nu)$ is obvious. Using the definition of $h^*$, $\widetilde{Q}_0$, and Theorems 3.1 and 3.2, we have

$$\Delta(h^*) = \tfrac{1}{2} \int_\Theta \int_\Theta \int_{\mathcal{X}} (h^*(\theta) - h^*(\eta))^2 Q(d\theta|x)Q(d\eta|x)M(dx)$$
$$= \tfrac{1}{2} \int_\Theta \int_\Theta \int_{\mathcal{X}} (\widetilde{h}(t(\theta)) - \widetilde{h}(t(\eta)))^2 Q(d\theta|x)Q(d\eta|x)M(dx)$$
$$= \tfrac{1}{2} \int_0^\infty \int_0^\infty \int_{\mathcal{X}} (\widetilde{h}(a) - \widetilde{h}(b))^2 \widetilde{Q}_0(da|x)\widetilde{Q}_0(db|x)M(dx) = \widetilde{\Delta}(\widetilde{h}). \quad \square$$

THEOREM 3.4. *If the $\widetilde{P}$-$s$ chain is l-$s$-r, then the $P$-$\nu$ chain is l-$\nu$-r.*

PROOF. Let $A$ be an $s$-proper subset of $[0, \infty)$. By Theorem 2.4,

$$\inf_{\widetilde{h} \in \widetilde{V}(A)} \widetilde{\Delta}(\widetilde{h}) = 0,$$

where $\widetilde{V}(\cdot)$ is the analog of $V(\cdot)$ [defined by (9)] for the $\widetilde{P}$-$s$ problem. Now, let $\{D_i\}$ be an increasing sequence of $s$-proper subsets of $[0, \infty)$ such that $D_i \nearrow [0, \infty)$. Setting $E_i = t^{-1}(D_i)$ for $i = 1, 2, \ldots$, we see that $\nu(E_i) = s(D_i) \in (0, \infty)$ and $E_i \nearrow \Theta$. By Theorem 2.2, it is sufficient to show that each $E_i$ is l-$\nu$-r. But using Theorem 2.4 shows that it is sufficient to verify that

(15) $$\inf_{h \in V(E_i)} \Delta(h) = 0.$$

To show (15), let $\varepsilon > 0$ be given and select $\widetilde{h} \in \widetilde{V}(D_i)$ so that $\widetilde{\Delta}(\widetilde{h}) < \varepsilon$. Then set $h^*(\theta) = \widetilde{h}(t(\theta))$. Clearly, $h^* \in V(E_i)$ and by Theorem 3.3, $\Delta(h^*) = \widetilde{\Delta}(\widetilde{h}) < \varepsilon$. Thus, (15) holds. □

REMARK 3.2. The choice of $[0, \infty)$ for the range of $t$ is purely for simplicity and the direct application to Example 1.1. Any other choice of Polish space for the range of $t$ yields the same results as long as Assumption 3.1 holds. A convenient choice for $t$ and its range is dependent upon the application.

**4. Recurrence of Markov chains.** In this section we discuss the recurrence of a discrete time Markov chain with values in a Polish space $\mathcal{X}$. The $\sigma$-algebra generated by the open sets of $\mathcal{X}$ is denoted by $\mathcal{B}$. The notation and setting in this section are independent of that in the earlier sections. Indeed, the material here may be of independent interest since some nontrivial generalizations of results in [17] are given below.

Consider a discrete time Markov chain $X = (X_0, X_1, \ldots)$ with state space $\mathcal{X}$ so each $X_i$ is an element of $\mathcal{X}$, $i = 0, 1, \ldots$. Therefore, $X$ takes values in the product space $\mathcal{X}^\infty$ which is equipped with the natural product $\sigma$-algebra $\mathcal{B}^\infty$. The one-step transition function of the chain, assumed to be a Markov kernel on $\mathcal{B} \times \mathcal{X}$, is denoted by $T(\cdot|y)$ for $y \in \mathcal{X}$. Given $x_0 \in \mathcal{X}$, let $P(\cdot|x_0)$ denote the distribution of $X \in \mathcal{X}^\infty$ when $X_0 = x_0$. Thus, $P(\cdot|x_0)$ is a probability measure on $\mathcal{B}^\infty$.

Now, let $C$ be a Borel subset of $\mathcal{X}$. Define the stopping time

$$\tau_C = \begin{cases} \infty, & \text{if } X_n \notin C \text{ for all } n \geq 1, \\ \text{smallest } n \geq 1 \text{ such that } X_n \in C, & \text{otherwise} \end{cases}$$

and let $E_C = \{w \in \mathcal{X}^\infty | \tau_C(w) < \infty\}$.

DEFINITION 4.1. The set $C$ is $P(\cdot|x_0)$-*recurrent* if $P(E_C|x_0) = 1$. The set $C$ is *recurrent* if $C$ is $P(\cdot|x_0)$-recurrent for all $x_0 \in \mathcal{X}$.

REMARK 4.1. If $C$ is recurrent and $\nu$-proper, in the language of Section 2, then $C$ is $\nu$-recurrent.

PROPOSITION 4.1. *If $C$ is $P(\cdot|x_0)$-recurrent for all $x_0 \in C^c$, then $C$ is recurrent.*

PROOF. This is a standard Markov chain argument which is omitted. □

4.1. *A condition for recurrence.* Let $f$ be a Borel measurable function defined on $\mathcal{X}$ with values in $[0, \infty)$. Define a sequence of random variables $Y_0, Y_1, Y_2, \ldots$ by

$$Y_n = f(X_{\tau_C \wedge n}), \qquad n = 0, 1, \ldots,$$

where $\wedge$ denotes the minimum. Obviously, $Y_i = f(X_i)$ for $i = 0, 1$. The proof of the following is omitted.

LEMMA 4.1. *With $E_C$ as defined above, if $w \in E_C$, then the sequence $\{Y_n(w) \mid n = 0, 1, 2, \ldots\}$ converges to a finite limit, namely, $f(X_{\tau_C(w)}(w))$.*

Here is an important structural condition on $f$ and the chain $X$.

ASSUMPTION 4.1. *Given $x_0 \in \mathcal{X}$,*

$$P\left(\limsup_{n \to \infty} f(X_n) = \infty \,\Big|\, x_0\right) = 1.$$

THEOREM 4.1. *Let $H \subseteq \mathcal{X}^\infty$ be the set of $w$'s such that $Y_n(w)$ converges to a finite limit. Suppose Assumption 4.1 holds and that $P(H|x_0) = 1$. Then the set $C$ is $P(\cdot \mid x_0)$-recurrent.*

PROOF. Let $E_C$ be as defined above. Then by Lemma 4.1, $E_C \subseteq H$. Thus $H = E_C \cup (H \cap E_C^c)$. To show that $C$ is $P(\cdot \mid x_0)$-recurrent, it suffices to show $P(H \cap E_C^c \mid x_0) = 0$ since then, $P(H|x_0) = P(E_C|x_0) = 1$. Let

$$K = \left\{ w \,\Big|\, \limsup_{n \to \infty} f(X_n(w)) = \infty \right\}.$$

Then, since $P(K|x_0) = 1$, we have $P(H \cap E_C^c \mid x_0) = P(H \cap E_C^c \cap K \mid x_0)$. Since $H \cap E_C^c \cap K$ is empty, this completes the proof. $\square$

The above shows that when Assumption 4.1 holds, $P(\cdot|x_0)$-recurrence of $C$ will hold if we can show that $\{Y_n | n = 1, \ldots\}$ converges to a finite limit a.s.-$P(\cdot|x_0)$. Of course, if $\{Y_n\}$ is a $P(\cdot|x_0)$ supermartingale, then $H$ has $P(\cdot|x_0)$ probability one and Theorem 4.1 applies. It is this martingale argument that [17] used.

4.2. *When is $\{Y_n\}$ a supermartingale?* Given the Markov chain $X = (X_0, X_1, \ldots)$ with state space $\mathcal{X}$ and $X_0 = x_0$, let $\mathcal{F}_n$ be the $\sigma$-algebra in $\mathcal{B}^\infty$ generated by $X_0, \ldots, X_n$. Recall that $\{Y_n, \mathcal{F}_n \mid n = 0, 1, 2 \ldots\}$ is a $P(\cdot|x_0)$ supermartingale if

(16) $$E(Y_{n+1}|\mathcal{F}_n) \leq Y_n, \qquad n = 0, 1, \ldots,$$

where the expectation is taken under $P(\cdot|x_0)$ on $\mathcal{B}^\infty$.

DEFINITION 4.2. A function $f:\mathcal{X} \to [0,\infty)$ is *superharmonic on $C^c$* if

(17) $$E(f(X_1) \mid X_0 = x_0) \leq f(x_0) \qquad \text{for all } x_0 \in C^c.$$

THEOREM 4.2. *If $f$ is superharmonic on $C^c$, then $(Y_n, \mathcal{F}_n \mid n = 0, 1, \ldots)$ is a supermartingale for each $x_0 \in C^c$.*

PROOF. Recall that $Y_i = f(X_i)$ for $i = 0, 1$ and fix $x_0 \notin C$. That (16) holds for $n = 0$ is a direct consequence of (17) since $x_0 \notin C$.

For $n \geq 1$, let $G_n \in \mathcal{F}_n$ be the event $G_n = \{X_1 \notin C, \ldots, X_n \notin C\} = \{\tau_C > n\}$. Obviously, $E(Y_{n+1} \mid \mathcal{F}_n) = E(I_{G_n^c} Y_{n+1} \mid \mathcal{F}_n) + E(I_{G_n} Y_{n+1} \mid \mathcal{F}_n)$. On the set $G_n^c$, $Y_{n+1} = Y_n$ so

(18) $$E(I_{G_n^c} Y_{n+1} \mid \mathcal{F}_n) = E(I_{G_n^c} Y_n \mid \mathcal{F}_n) = I_{G_n^c} Y_n.$$

On the set $G_n$, $\tau_C \geq n + 1$ so $Y_{n+1} = f(X_{n+1})$ and

$$E(I_{G_n} Y_{n+1} \mid \mathcal{F}_n) = E(I_{G_n} f(X_{n+1}) \mid \mathcal{F}_n) = I_{G_n} E(f(X_{n+1}) \mid \mathcal{F}_n) \leq I_{G_n} f(X_n).$$

The last inequality follows from the Markov property and the assumption that $f$ is superharmonic on $C^c$. Thus, on the set $G_n$,

(19) $$E(I_{G_n} Y_{n+1} \mid \mathcal{F}_n) \leq I_{G_n} f(X_n) = I_{G_n} Y_n.$$

Combining (18) and (19) shows (16) holds. □

COROLLARY 4.1. *Suppose $f$ is superharmonic on $C^c$. Then for each $x_0 \in C^c$, $\{Y_n\}$ converges a.s.-$P(\cdot \mid x_0)$ to a finite random variable $Y$.*

PROOF. This is just the Supermartingale Convergence Theorem since $0 \leq EY_n \leq EY_1$ for all $n$, so $0 \leq \sup_n EY_n \leq EY_1$ ([3], page 468). □

Combining what has now been established, we have the following result.

THEOREM 4.3. *Suppose Assumption 4.1 holds for $x_0 \in C^c$ and that $f$ is superharmonic on the set $C^c \subseteq \mathcal{X}$. Then the set $C$ is recurrent.*

Theorem 4.3 parallels some results in Chapter 8 of [19] in the use of superharmonic functions. An important difference is that Meyn and Tweedie assume the Markov chain is irreducible, while we rely on Assumption 4.1.

4.3. *The case when* $\mathcal{X} = [0, \infty)$  In this section it is assumed that the state space $\mathcal{X}$ is $[0, \infty)$. We develop sufficient conditions on the chain $X = (X_0, X_1, \ldots)$ so that a particular function $f$ is superharmonic on the set $[m, \infty)$ when $m$ is large enough. In addition, we provide a sufficient condition for Assumption 4.1 which may be easy to check in examples. When these two results hold, it will follow from our previous results that the set $C = [0, m)$ is recurrent. An application of this result is provided in Section 5.

We begin with a statement of a main result, although several auxiliary results are needed before the proof can be completed. For $k = 1, 2, 3$, let

$$\mu_k(x) = \int_0^\infty (y-x)^k T(dy|x). \tag{20}$$

The quantity

$$\int_0^\infty |y-x|^3 T(dy|x)$$

is assumed finite for all $x \in [0, \infty)$, so $\mu_k(x)$ is well defined for all $x \in [0, \infty)$ and $k = 1, 2, 3$. We also assume that $\mu_2(x) > 0$ for all sufficiently large $x$.

THEOREM 4.4. *Assume that there is a function $\psi_1$ such that, for all sufficiently large $x$,*

$$\mu_1(x) \leq \frac{\mu_2(x)}{2x}[1 + \psi_1(x)] \tag{21}$$

*and*

$$\lim_{x \to \infty} (\log x)\psi_1(x) = 0. \tag{22}$$

*Also, assume that*

$$\lim_{x \to \infty} \frac{\log x}{x} \frac{\mu_3(x)}{\mu_2(x)} = 0. \tag{23}$$

*Then, for $x \geq 0$,*

$$f_0(x) = \log(\log(e + x)) \tag{24}$$

*is superharmonic on the interval $[m, \infty)$ for $m$ large enough.*

What needs to be established to prove Theorem 4.4 is that, for all large $x$,

$$\delta(x) = E(f_0(X_1) - f_0(x) \mid X_0 = x) \leq 0. \tag{25}$$

To this end, note that the first four derivatives of $f_0$ satisfy

$$f_0'(x) = \frac{1}{(x+e)\log(x+e)} > 0,$$

$$f_0''(x) = -\frac{\log(x+e)+1}{(x+e)^2\log^2(x+e)} < 0,$$

$$f_0'''(x) = \frac{2\log^2(x+e) + 3\log(x+e) + 2}{(x+e)^3\log^3(x+e)} > 0,$$

$$f_0^{(iv)}(x) = -\frac{6\log^3(x+e) + 7\log^2(x+e) + 12\log(x+e) + 6}{(x+e)^4\log^4(x+e)} < 0.$$

Now, expanding $f_0$ in a Taylor series about $x$ in (25), discarding the negative term $f^{(iv)}$ in this expansion, and doing a bit of algebra results in

$$(26) \qquad \delta(x) \leq \frac{f_0'(x)\mu_2(x)}{2x}\left[\frac{2x\mu_1(x)}{\mu_2(x)} + \frac{f_0''(x)}{f_0'(x)}x + \frac{2f_0'''(x)\mu_3(x)x}{6f_0'(x)\mu_2(x)}\right].$$

For notational convenience, set

$$\psi_2(x) = \frac{2f_0'''(x)\mu_3(x)x}{6f'(x)\mu_2(x)}.$$

Our first intermediate conclusion is the following, which is a direct consequence of (26) and (21).

LEMMA 4.2. *For all sufficiently large $x$, with $\delta(x)$ as defined in (25),*

$$(27) \qquad \delta(x) \leq \frac{f_0'(x)\mu_2(x)}{2x}\left[1 + \psi_1(x) + \frac{f_0''(x)x}{f_0'(x)} + \psi_2(x)\right].$$

LEMMA 4.3. *The right-hand side of the inequality (27) is equal to*

$$(28) \qquad \frac{f_0'(x)\mu_2(x)}{2x}xf_0'(x)\left[-1 + \frac{e\log(x+e)}{x} + \frac{\psi_1(x)}{xf_0'(x)} + \frac{\psi_2(x)}{xf_0'(x)}\right].$$

PROOF. First verify by direct calculation that

$$(29) \qquad \frac{f_0'(x) + xf_0''(x)}{(f_0'(x))^2} = -x + e\log(x+e).$$

Multiplying and dividing the bracketed term in (27) by $xf_0'(x)$ and using (29) immediately yields the claim. □

Now, we complete the proof of Theorem 4.4 by showing that the bracketed term in (28) is negative for all large $x$. But, this is a direct consequence of assumptions (22) and (23). For example,

$$\frac{\psi_1(x)}{xf_0'(x)} = \frac{x+e}{x}\log(x+e)\psi_1(x),$$

which obviously converges to zero as $x \to \infty$ under assumption (22). With a bit more algebra, (23) implies that $\psi_2(x)/xf_0'(x)$ converges to zero as $x \to \infty$. Therefore, the right-hand side of (27) is negative for all large $x$. Thus, $\delta(x) \leq 0$ for all $x$ large enough and the proof is complete.

We now turn to a discussion of the structural Assumption 4.1 when $\mathcal{X} = [0, \infty)$ and the function $f$ is $f_0$ in (24). Because $f_0$ is monotone increasing from $[0, \infty)$ onto $[0, \infty)$, it is clear that Assumption 4.1 holds for $f = f_0$ if and only if the following holds:

CONDITION 4.1. Given $x_0 \in [0, \infty)$,
$$P\left(\limsup_{n \to \infty} X_n = \infty \mid x_0\right) = 1.$$

PROPOSITION 4.2. *Define the event*
$$A_{m,k} = \{X_1 \in [0, m], X_2 \in [0, m], \ldots, X_k \in [0, m]\}$$
*for positive integers $m$ and $k$. Assume that, for each $x_0$ and each $m \in \mathbb{N}$, $P(A_{m,k}|x_0) \to 0$ as $k \to \infty$. Then Condition 4.1 holds for all $x_0$.*

PROOF. Fix $x_0$. For $m \in \mathbb{N}$, define $E_m = \{X_n \in [0, m] \text{ for } n = 1, 2, \ldots\}$. Let $E = \bigcup_{m=1}^{\infty} E_m$. Clearly, $P(E|x_0) = 0$ if and only if $P(E_m|x_0) = 0$ for all $m \in \mathbb{N}$. Define
$$F = \left\{\limsup_{n \to \infty} X_n = \infty\right\},$$
and note that $E^c = F$. Hence, $P(F|x_0) = 1$ if and only $P(E_m|x_0) = 0$ for all $m \in \mathbb{N}$. Note that $A_{m,k} \downarrow E_m$. By assumption, for each fixed $m \in \mathbb{N}$, $P(A_{m,k}|x_0) \to 0$ as $k \to \infty$. Thus, $P(E_m|x_0) = 0$ for all $m \in \mathbb{N}$ and the result is proved. □

PROPOSITION 4.3. *Assume that for each positive integer $m$ there exists a $\delta = \delta(m) < 1$ such that*
$$(30) \qquad \sup_{x_0 \in [0,m]} T([0, m] \mid x_0) \leq \delta.$$
*Then Condition 4.1 holds for all $x_0$.*

PROOF. Fix $x_0$ and $m \in \mathbb{N}$. Note that
$$P(A_{m,k}|x_0) = \int_0^m \int_0^m \cdots \int_0^m T(dx_1|x_0) \cdots T(dx_{k-1}|x_{k-2}) T([0,m]|x_{k-1})$$
$$\leq T([0, m] \mid x_0) \delta^{k-1}$$
$$\leq \delta^{k-1}.$$

Thus, $P(A_{m,k}|x_0) \to 0$ as $k \to \infty$ and an application of Proposition 4.2 completes the proof. □

In summary, the main conclusion of this section is the following.

THEOREM 4.5. *Assume that $\mathcal{X} = [0, \infty)$ and that the assumptions of Theorem 4.4 and Proposition 4.3 hold. Then there exists an $m \in (0, \infty)$ such that the set $C = [0, m)$ is recurrent.*

It is Theorem 4.5 that is used in the application of the next section.

**5. Strongly admissible priors for the multivariate normal mean.** We now use our results to identify strongly admissible priors for the mean of a multivariate normal distribution. Recall the setting of Example 1.1. Assume that $X \sim N_p(\theta, I_p)$ and take the prior to be $\nu_{a,b}(d\theta) = d\theta/(a + \|\theta\|^2)^b$, where $a \geq 0$, $b > 0$, and $d\theta$ is Lebesgue measure on $\mathbb{R}^p$.

The prior $\nu_{0,b}$ is improper for all $b > 0$, but the marginal is $\sigma$-finite only when $b < p/2$. On the other hand, when $a > 0$, $\nu_{a,b}$ is improper only when $b \leq p/2$ and the marginal is $\sigma$-finite for all $b$ in this range. We therefore restrict attention to $b \in (0, p/2)$ when $a = 0$ and to $b \in (0, p/2]$ when $a > 0$. When these conditions are satisfied, a proper posterior distribution $Q_{a,b}(d\theta|x)$ exists and satisfies the disintegration

$$P(dx|\theta)\nu_{a,b}(d\theta) = Q_{a,b}(d\theta|x)M_{a,b}(dx).$$

Here is the main result of this section.

THEOREM 5.1. *Suppose $X \sim N_p(\theta, I_p)$ with $p \geq 3$. The prior $\nu_{a,b}$ [equivalently, the posterior $Q_{a,b}(d\theta|x)$] is strongly admissible if either (A) $a > 0$ and $b \in [p/2 - 1, p/2]$ or (B) $a = 0$ and $b \in [p/2 - 1, p/2)$.*

REMARK 5.1. The problem we address is substantively different from the traditional problem of estimating $\theta$ under quadratic loss. However, one of the results established in [2] is that when $a = 1$ and $b = (p-1)/2$ the formal Bayes estimator of $\theta$ is admissible under quadratic loss. The argument in [2] uses the general results in [4]. It seems plausible that the results in [4] could be used to obtain admissibility, under quadratic loss, of the formal Bayes estimator of $\theta$ for any of the priors $\nu_{a,b}$ in Theorem 5.1.

PROOF OF THEOREM 5.1. Throughout, we will refrain from using subscripts on $\nu$, $M$ and $Q$. In Example 3.1 it is shown that the prior $\nu(d\theta)$ can be expressed as $\pi(d\theta|\beta)s(d\beta)$, where $\pi(\cdot|\beta)$ is the uniform distribution on $\{\theta|\|\theta\|^2 = \beta\}$ and $s(d\beta)$ is given in (12) with $g_0(\beta) = (a+\beta)^{-b}$. The new model, $\widetilde{P}$, has density (with respect to Lebesgue measure on $\mathbb{R}^p$) given by

$$\tilde{f}(x|\beta) = \int_\Xi (2\pi)^{-p/2} e^{-(1/2)\|x-\sqrt{\beta}\xi\|^2} \pi_1(d\xi),$$

where $\pi_1$ is the uniform distribution on $\Xi = \{\theta|\|\theta\| = 1\}$. The formal posterior for $\beta$ can be written as $\widetilde{Q}(d\beta|x) = q(\beta|x)\,d\beta$, where

$$(31) \qquad q(\beta|x) = \frac{c\tilde{f}(x|\beta)g_0(\beta)\beta^{p/2-1}}{m(x)}$$

with $c$ a positive constant and

$$(32) \qquad m(x) = c\int_0^\infty \tilde{f}(x|\beta)g_0(\beta)\beta^{p/2-1}\,d\beta = \int_{\mathbb{R}^p} f(x|\theta)g_0(\|\theta\|^2)\,d\theta,$$

where $f(x|\theta)$ is the multivariate normal density with mean $\theta \in \mathbb{R}^p$ and covariance matrix the $p \times p$ identity, $I_p$. The $\widetilde{P}$-$s$ chain has Markov transition function given by

$$(33) \qquad \widetilde{R}(d\beta|\eta) = \int_{\mathbb{R}^p} \widetilde{Q}(d\beta|x)\tilde{f}(x|\eta)\,dx.$$

According to Theorem 3.4, to prove that $\nu$ is strongly admissible, it suffices to show that the $\widetilde{P}$-$s$ chain is l-$s$-r. This is now established by showing that the conditions of Theorem 4.5 are satisfied.

Note that $\tilde{f}(x|\eta)$ and $q(\beta|x)$ are strictly positive for all $\beta, \eta \in (0,\infty)$ and $x \in \mathbb{R}^p$. It follows that $\widetilde{R}(C|\eta) > 0$ for every $\eta \in [0,\infty)$ and every $C$ with positive Lebesgue measure. In order to apply Proposition 4.3, we need to show that, for any $m \in \mathbb{N}$, there exists a $\delta < 1$ such that

$$\sup_{\eta \in [0,m]} \widetilde{R}([0,m]|\eta) \leq \delta.$$

It suffices to show that $\widetilde{R}([0,m]|\eta)$ is a continuous function of $\eta$ for $\eta \in [0,m]$. Fix $\eta^* \in [0,m]$ and, without loss of generality, let $\{\eta_k\}_{k=1}^\infty$ be a sequence in $[0,4m]$ that converges to $\eta^*$. An application of dominated convergence shows that $\tilde{f}(x|\eta_k) \to \tilde{f}(x|\eta^*)$ as $k \to \infty$. Let $g\colon \mathbb{R}^p \to \mathbb{R}$ such that

$$g(x) = I(\|x\| < 4\sqrt{m}) + I(\|x\| \geq 4\sqrt{m})e^{-(1/2)\|x/2\|^2}.$$

Since $g$ is integrable and $\tilde{f}(x|\eta_k) \leq g(x)$ for all $x$ and all $k$, another application of dominated convergence yields the desired continuity.

We now turn our attention to establishing the conditions of Theorem 4.4. Recall that

$$\mu_k(\eta) := \int_0^\infty (\beta - \eta)^k \widetilde{R}(d\beta|\eta).$$

In the appendix we prove

$$\mu_1(\eta) = 2p - 4b + \psi_1^*(\eta), \qquad \mu_2(\eta) = 8\eta + \psi_2^*(\eta) \quad \text{and} \quad \mu_3(\eta) = \psi_3^*(\eta),$$

where, as $\eta \to \infty$, $\psi_1^*(\eta) = O(\eta^{-1})$, $\psi_2^*(\eta) = O(1)$ and $\psi_3^*(\eta) = O(\eta)$. Set $\psi_1(\eta) = \eta^{-\varepsilon}$ for some $\varepsilon \in (0,1)$. Then (22) holds and

$$\mu_2(\eta)[1 + \psi_1(\eta)] = 8\eta + 8\eta^{1-\varepsilon} + O(1).$$

Since $b \geq p/2 - 1$, it follows that

$$2\eta\mu_1(\eta) = (4p - 8b)\eta + O(1) \leq 8\eta + O(1),$$

so (21) holds. Furthermore, (23) holds since, as $\eta \to \infty$,

$$\frac{\log \eta}{\eta} \frac{\mu_3(\eta)}{\mu_2(\eta)} = \frac{\log \eta}{\eta} \frac{\psi_3^*(\eta)}{8\eta + \psi_2^*(\eta)} = \frac{\psi_3^*(\eta)}{\eta} \frac{\log \eta}{8\eta + \psi_2^*(\eta)} \to 0. \qquad \square$$

Let $\nu(d\theta)$ be an improper prior on $\Theta = \mathbb{R}^p$ and suppose that the $P$-$\nu$ chain is l-$\nu$-r so that $\nu$ is strongly admissible. A result in [11] shows that (under mild conditions) if $h: \mathbb{R}^p \to [0, \infty)$ is a bounded function, then the "perturbed" prior $\nu^*(d\theta) = h(\theta)\nu(d\theta)$ is also strongly admissible. In fact, Corollary 4 in [11] in conjunction with the results in the proof of our Theorem 5.1 immediately yields the following:

THEOREM 5.2. *Suppose $X \sim N_p(\theta, I_p)$ with $p \geq 3$. Let $h: \mathbb{R}^p \to [0, \infty)$ be bounded and suppose the perturbed prior $\nu_{h,a,b}(d\theta) = h(\theta)\nu_{a,b}(d\theta)$ is improper. Then $\nu_{h,a,b}$ is strongly admissible if either (A) $a > 0$ and $b \in [p/2 - 1, p/2]$ or (B) $a = 0$ and $b \in [p/2 - 1, p/2)$.*

REMARK 5.2. It is possible to get most of part (A) of Theorem 5.1 by combining the proof of part (B) with the work in [11]. Fix $b \in [p/2 - 1, p/2)$. We know from the proof of part (B) that the Markov chain associated with the prior $\nu_{0,b}(d\theta) = d\theta/(\|\theta\|^2)^b$ is l-$\nu_{0,b}$-r. Now fix $a > 0$, and note that

$$\nu_{a,b}(d\theta) = \frac{1}{(a + \|\theta\|^2)^b} d\theta = \left(\frac{\|\theta\|^2}{a + \|\theta\|^2}\right)^b \nu_{0,b}(d\theta).$$

Since $(\|\theta\|^2/(a + \|\theta\|^2))^b$ is a bounded function, the results in [11] imply that $\nu_{a,b}$ is strongly admissible.

Finally, using arguments somewhat similar to those in the proof of Theorem 5.1, the multivariate Poisson case was discussed in [16].

## APPENDIX: ASYMPTOTIC EXPANSION OF MOMENTS.

Throughout this section, $p \in \{3, 4, \ldots\}$, $a \geq 0$ and $b > 0$ are considered fixed and the assumptions of Theorem 5.1 are in force. We begin by recalling some facts concerning the noncentral $\chi^2$ distribution. Suppose that $Z \sim N_p(\gamma, I_p)$ and let $\lambda = \|\gamma\|^2$. Then $V = \|Z\|^2 \sim \chi_p^2(\lambda)$ and its density at $v > 0$ is given by

$$\text{(34)} \qquad \sum_{n=0}^{\infty} \frac{e^{-\lambda/2}(\lambda/2)^n}{n!} \frac{(v/2)^{n+p/2-1} e^{-v/2}}{2\Gamma(n+p/2)}.$$

We will need the first three moments of $V$, which are as follows:

$$E[V|\lambda] = \lambda + p, \qquad E[V^2|\lambda] = (\lambda + p)^2 + 4\lambda + 2p,$$
$$E[V^3|\lambda] = (\lambda + p)^3 + 12(\lambda + p)^2 - 6\lambda p + 24\lambda - 6p^2 + 8p.$$

Recall that $g_0(z) = (a+z)^{-b}$. Using (34), one can show that $E[g_0(V)V^k|\lambda] = 2^k E[w_k(N)|\lambda]$, where $N|\lambda \sim \text{Poisson}(\lambda/2)$ and, for $n \in \mathbb{Z}^+ := \{0, 1, 2, \ldots\}$,

$$w_k(n) := \int_0^\infty \frac{g_0(z)(z/2)^{n+p/2+k-1} e^{-z/2}}{2\Gamma(n+p/2)} \, dz.$$

When $a = 0$ and $p/2 - b > 0$,

$$\text{(35)} \qquad w_k(n) = \frac{\Gamma(n+p/2+k-b)}{2^b \Gamma(n+p/2)},$$

which is well defined even when $n = k = 0$. Our goal is to prove the following.

PROPOSITION A.1. *For $\eta > 0$, we have*

$$\mu_1(\eta) = 2p - 4b + \psi_1^*(\eta), \qquad \mu_2(\eta) = 8\eta + \psi_2^*(\eta) \quad \text{and} \quad \mu_3(\eta) = \psi_3^*(\eta),$$

*where, as $\eta \to \infty$, $\psi_1^*(\eta) = O(\eta^{-1})$, $\psi_2^*(\eta) = O(1)$ and $\psi_3^*(\eta) = O(\eta)$.*

PROOF. The proof is constructed in several intermediate steps. Define $m_k(\zeta) \colon \mathbb{R}_+ \to \mathbb{R}_+$ via $m_k(\zeta) = E[g_0(U)U^k|\zeta]$ with $U \sim \chi_p^2(\zeta)$. Our first result follows.

PROPOSITION A.2. *Fix $\eta > 0$ and let $Y \sim \chi_p^2(\eta)$. For $k \in \{0, 1, 2, 3\}$,*

$$\text{(36)} \qquad E\left[\frac{m_k(Y)}{m_0(Y)}\bigg|\eta\right] < \infty$$

*and*

$$\text{(37)} \qquad \int_0^\infty \beta^k \widetilde{R}(d\beta|\eta) = E\left[\frac{m_k(Y)}{m_0(Y)}\bigg|\eta\right].$$

PROOF. We begin with (36). Since $g_0$ is convex, we can use Jensen's inequality to obtain

$$\frac{m_k(\zeta)}{m_0(\zeta)} \leq \frac{m_k(\zeta)}{g_0(\zeta+p)} = (a+\zeta+p)^b m_k(\zeta).$$

Hence, it suffices to show that $E[Y^b m_k(Y) \mid \eta] < \infty$, but this follows from the fact that $Y$ has a moment generating function.

We now establish (37). Recall that $\pi_1$ is the uniform distribution on $\Xi = \{\theta \mid \|\theta\| = 1\}$. Then

$$\int_0^\infty \beta^k \tilde{f}(x|\beta) g_0(\beta) \beta^{p/2-1} \, d\beta$$

$$= \int_0^\infty \beta^k \int_\Xi (2\pi)^{-p/2} e^{-(1/2)\|x - \sqrt{\beta}\xi\|^2} \pi_1(d\xi) g_0(\beta) \beta^{p/2-1} \, d\beta$$

$$= \int_{\mathbb{R}^p} g_0(\|\theta\|^2)(\|\theta\|^2)^k (2\pi)^{-p/2} e^{-(1/2)\|x-\theta\|^2} \, d\theta.$$

The last expression is equal to

$$E[g_0(\|\theta\|^2)(\|\theta\|^2)^k],$$

where $\theta \sim \mathrm{N}_p(x, I_p)$. Of course, if $\theta \sim \mathrm{N}_p(x, I_p)$, then $\|\theta\|^2 \sim \chi_p^2(\|x\|^2)$. Thus,

$$\int_0^\infty \beta^k \tilde{f}(x|\beta) g_0(\beta) \beta^{p/2-1} \, d\beta = m_k(\|x\|^2).$$

Now note that

$$\int_0^\infty \beta^k \widetilde{Q}(d\beta|x) = \frac{\int_0^\infty \beta^k \tilde{f}(x|\beta) g_0(\beta) \beta^{p/2-1} d\beta}{\int_0^\infty \tilde{f}(x|\beta) g_0(\beta) \beta^{p/2-1} d\beta} = \frac{m_k(\|x\|^2)}{m_0(\|x\|^2)}.$$

To complete the proof, note that

$$\int_0^\infty \beta^k \widetilde{R}(d\beta|\eta) = \int_{\mathbb{R}^p} \left[ \int_0^\infty \beta^k \widetilde{Q}(d\beta|x) \right] \tilde{f}(x|\eta) \, dx = \int_{\mathbb{R}^p} \frac{m_k(\|x\|^2)}{m_0(\|x\|^2)} \tilde{f}(x|\eta) \, dx$$

$$= \int_{\mathbb{R}^p} \int_\Xi \frac{m_k(\|x\|^2)}{m_0(\|x\|^2)} \left(\frac{1}{\sqrt{2\pi}}\right)^p e^{-(1/2)\|x-\sqrt{\eta}\xi\|^2} \pi_1(d\xi) \, dx$$

$$= \int_\Xi \left[ \int_{\mathbb{R}^p} \frac{m_k(\|x\|^2)}{m_0(\|x\|^2)} \left(\frac{1}{\sqrt{2\pi}}\right)^p e^{-(1/2)\|x-\sqrt{\eta}\xi\|^2} \, dx \right] \pi_1(d\xi).$$

Consider the inside integral. If $X \sim \mathrm{N}_p(\sqrt{\eta}\xi, I_p)$, then, since $\|\sqrt{\eta}\xi\|^2 = \eta$, it follows that $\|X\|^2 \sim \chi_p^2(\eta)$. Hence, the inside integral can be expressed as

$$E\left[\frac{m_k(Y)}{m_0(Y)} \bigg| \eta\right],$$

where $Y \sim \chi_p^2(\eta)$, and this does not depend on $\xi$. Therefore,

$$\int_0^\infty \beta^k \widetilde{R}(d\beta|\eta) = \int_\Xi E\left[\frac{m_k(Y)}{m_0(Y)}\bigg|\eta\right]\pi(d\xi) = E\left[\frac{m_k(Y)}{m_0(Y)}\bigg|\eta\right]. \qquad \square$$

It now follows from the discussion at the beginning of this appendix that $m_k(y) = 2^k E[w_k(N)|y]$, where $N|y \sim \text{Poisson}(y/2)$ and

$$\mu_1(\eta) = 2E\left[\frac{E(w_1(N)|Y)}{E(w_0(N)|Y)}\bigg|\eta\right] - \eta,$$

$$\mu_2(\eta) = 4E\left[\frac{E(w_2(N)|Y)}{E(w_0(N)|Y)}\bigg|\eta\right] - 4\eta E\left[\frac{E(w_1(N)|Y)}{E(w_0(N)|Y)}\bigg|\eta\right] + \eta^2,$$

$$\mu_3(\eta) = 8E\left[\frac{E(w_3(N)|Y)}{E(w_0(N)|Y)}\bigg|\eta\right] - 12\eta E\left[\frac{E(w_2(N)|Y)}{E(w_0(N)|Y)}\bigg|\eta\right]$$

$$+ 6\eta^2 E\left[\frac{E(w_1(N)|Y)}{E(w_0(N)|Y)}\bigg|\eta\right] - \eta^3,$$

where $N|Y = y \sim \text{Poisson}(y/2)$ and $Y \sim \chi_p^2(\eta)$.

We now begin working on $w_k(n)$. Our first result is as follows.

PROPOSITION A.3. *Suppose $n \in \mathbb{Z}^+$ and $k \in \{1,2,3\}$.*

1. *If $a = 0$, then*

$$\frac{w_k(n)}{w_{k-1}(n)} = n + \frac{p}{2} + k - b - 1. \tag{38}$$

2. *If $a > 0$, then*

$$\frac{w_k(n)}{w_{k-1}(n)} = n + \frac{p}{2} + k - b - 1 + \phi(n + k - 1), \tag{39}$$

*where $\phi : \mathbb{Z}^+ \to (0, \infty)$ is given by*

$$\phi(n) := \frac{\int_0^\infty (ab/(a+z))g_0(z)(z/2)^{n+p/2-1}e^{-z/2}\,dz}{\int_0^\infty g_0(z)(z/2)^{n+p/2-1}e^{-z/2}\,dz}.$$

*Moreover, $\phi(n)$ is bounded and $\phi(n) = O(n^{-1})$ as $n \to \infty$.*

PROOF. The $a = 0$ result follows directly from (35). Now assume that $a > 0$ and define $\kappa(z) = ab/(a+z)$. Note that

$$zg_0'(z) = [-b + \kappa(z)]g_0(z), \tag{40}$$

where $g_0'$ denotes the derivative of $g_0$. Integration by parts yields

$$w_k(n) = \frac{1}{\Gamma(n+p/2)} \int_0^\infty g_0'(z) e^{-z/2} \left(\frac{z}{2}\right)^{n+p/2+k-1} dz$$
$$+ \left(n + \frac{p}{2} + k - 1\right) w_{k-1}(n).$$

An application of (40) yields

$$\int_0^\infty g_0'(z) e^{-z/2} \left(\frac{z}{2}\right)^{n+p/2+k-1} dz$$
$$= \frac{1}{2} \int_0^\infty [-b + \kappa(z)] g_0(z) e^{-z/2} \left(\frac{z}{2}\right)^{n+p/2+k-2} dz$$
$$= -b\Gamma\left(n + \frac{p}{2}\right) w_{k-1}(n) + \frac{1}{2} \int_0^\infty \kappa(z) g_0(z) e^{-z/2} \left(\frac{z}{2}\right)^{n+p/2+k-2} dz.$$

It follows that

$$\frac{w_k(n)}{w_{k-1}(n)} = n + \frac{p}{2} + k - b - 1 + \phi(n + k - 1).$$

Since $\kappa(z)$ is bounded above by $b$, it is clear that $\phi(n)$ is also bounded above by $b$. Now define $N(n)$ and $D(n)$ as follows:

$$\phi(n) = \frac{N(n)}{D(n)} = \frac{\int_0^\infty (ab/(a+z)) g_0(z) (z/2)^{n+p/2-1} e^{-z/2} dz}{\int_0^\infty g_0(z) (z/2)^{n+p/2-1} e^{-z/2} dz}.$$

Using Jensen's inequality, we obtain

$$D(n) = \int_0^\infty g_0(z) \left(\frac{z}{2}\right)^{n+p/2-1} e^{-z/2} dz$$
$$= 2\Gamma(n+p/2) \int_0^\infty g_0(z) \frac{z^{n+p/2-1} e^{-z/2}}{\Gamma(n+p/2) 2^{n+p/2}} dz$$
$$\geq \frac{2\Gamma(n+p/2)}{(a+p+2n)^b}.$$

As long as $n > b + 1 - p/2$, we have

$$N(n) = \int_0^\infty \left(\frac{ab}{a+z}\right) g_0(z) \left(\frac{z}{2}\right)^{n+p/2-1} e^{-z/2} dz$$
$$\leq \frac{ab}{2^{n+p/2-1}} \int_0^\infty z^{n+p/2-b-2} e^{-z/2} dz$$
$$= \frac{ab\Gamma(n+p/2-b-1)}{2^b}.$$

Putting these two bounds together, we have

$$n\phi(n) \le \frac{ab}{2^{b+1}}\left(\frac{a+p}{n}+2\right)^b n^{b+1}\frac{\Gamma(n+p/2-b-1)}{\Gamma(n+p/2)}.$$

Finally, according to [1], page 257,

$$n^{b+1}\frac{\Gamma(n+p/2-b-1)}{\Gamma(n+p/2)} \to 1 \quad \text{as } n \to \infty,$$

from which it follows that $\phi(n) = O(n^{-1})$. □

A simple calculation shows that, for $n \in \mathbb{Z}^+$ and $k \in \{1,2,3\}$, we have

$$w_k(n) = \left(n+\frac{p}{2}\right)w_{k-1}(n+1).$$

This combined with (39) yields, when $a > 0$,

(41) $\quad (n+p/2)w_{k-1}(n+1) = (n+p/2+k-b-1+\phi(n+k-1))w_{k-1}(n).$

PROPOSITION A.4. *Suppose* $N|y \sim Poisson(y/2)$ *where* $y > 0$ *and let* $k \in \{1,2,3\}$.

1. *If* $a = 0$, *then*

(42)
$$\frac{E[w_k(N)|y]}{E[w_{k-1}(N)|y]}$$
$$= \frac{p}{2}+k-b-1+\frac{y}{2}+\frac{y(k-b-1)}{2}\frac{E[w_{k-1}(N)/(N+p/2)|y]}{E[w_{k-1}(N)|y]}.$$

2. *If* $a > 0$, *then*

(43)
$$\frac{E[w_k(N)|y]}{E[w_{k-1}(N)|y]}$$
$$= \frac{p}{2}+k-b-1+\frac{y}{2}+\frac{y(k-b-1)}{2}\frac{E[w_{k-1}(N)/(N+p/2)|y]}{E[w_{k-1}(N)|y]}$$
$$+\frac{y}{2}\frac{E[\phi(N+k-1)w_{k-1}(N)/(N+p/2)|y]}{E[w_{k-1}(N)|y]}$$
$$+\frac{E[\phi(N+k-1)w_{k-1}(N)|y]}{E[w_{k-1}(N)|y]}.$$

PROOF. First,

$$E[Nw_{k-1}(N)|y] = \sum_{n=0}^{\infty} nw_{k-1}(n)\frac{e^{-y/2}(y/2)^n}{n!}$$

$$= \frac{y}{2} \sum_{n=1}^{\infty} w_{k-1}(n) \frac{e^{-y/2}(y/2)^{n-1}}{(n-1)!}$$
$$= \frac{y}{2} \sum_{m=0}^{\infty} w_{k-1}(m+1) \frac{e^{-y/2}(y/2)^{m}}{m!}$$
$$= \frac{y}{2} E[w_{k-1}(N+1)|y].$$

We now prove the result for $a > 0$ and we note that a proof for the $a = 0$ case can be constructed simply by replacing $\phi$ by 0 in the following argument. Rearranging (39) and taking expectations yields

$$E[w_k(N)|y] = \left(\frac{p}{2} + k - b - 1\right) E[w_{k-1}(N)|y]$$
(44)
$$+ E[N w_{k-1}(N)|y] + E[\phi(N + k - 1) w_{k-1}(N)|y].$$

A rearrangement of (41) yields

$$w_{k-1}(n+1) = w_{k-1}(n) + w_{k-1}(n) \left[ \frac{2(k - b - 1 + \phi(n + k - 1))}{2n + p} \right].$$

Therefore,

$$E[N w_{k-1}(N)|y] = \frac{y}{2} E[w_{k-1}(N+1)|y]$$
(45)
$$= \frac{y}{2} E[w_{k-1}(N)|y]$$
$$+ y E\left[ \frac{(k - b - 1 + \phi(N + k - 1)) w_{k-1}(N)}{2N + p} \bigg| y \right].$$

Replacing $E[N w_{k-1}(N)|y]$ in (44) with the last line in (45) and dividing through by $E[w_{k-1}(N)|y]$ yields (43). $\square$

PROPOSITION A.5. *There exists a constant $c \in (0, \infty)$ such that*

$$\frac{w_{k-1}(n)}{w_{k-1}(n+1)} \leq c$$

*for all $n \in \mathbb{Z}^+$ and all $k \in \{1, 2, 3\}$. Note that $c$ may depend on $p$, $a$ and $b$, but does not depend on $n$ and $k$.*

PROOF. We first handle the case $a = 0$. Equation (35) shows that

$$\frac{w_{k-1}(n)}{w_{k-1}(n+1)} = \frac{n + p/2}{n + p/2 + k - b - 1}.$$

Since $b < p/2$, the denominator is always strictly positive. Furthermore, for any $k \in \{1, 2, 3\}$, this fraction clearly converges to 1 as $n \to \infty$. Thus, the fraction is bounded above by a positive, finite constant.

The $a > 0$ case is similar, but we do not have the luxury of using the exact expression (35). A rearrangement of (41) yields

$$\frac{w_{k-1}(n)}{w_{k-1}(n+1)} = \frac{n+p/2}{n+p/2+k-b-1+\phi(n+k-1)}.$$

Since $\phi(n)$ is strictly positive, the denominator is always strictly positive. As above, for any $k \in \{1,2,3\}$, the fraction converges to 1 as $n \to \infty$. Thus, the fraction is bounded above by a positive, finite constant. $\square$

Note that for all $n \in \mathbb{Z}^+$ and all $k \in \{1,2,3\}$, we have

$$\frac{w_{k-1}(n)}{w_{k-1}(n+2)} = \frac{w_{k-1}(n)}{w_{k-1}(n+1)} \frac{w_{k-1}(n+1)}{w_{k-1}(n+2)} \leq c^2.$$

PROPOSITION A.6. *Suppose $h:\mathbb{Z}^+ \to \mathbb{R}$ is such that $h(n) = O(n^{-l})$ as $n \to \infty$, where $l \in \{1,2\}$. If $N|y \sim \text{Poisson}(y/2)$ where $y > 0$, then there exists $0 < d < \infty$ such that for any $k \in \{1,2,3\}$,*

$$\left|\frac{E[h(N)w_{k-1}(N)|y]}{E[w_{k-1}(N)|y]}\right| \leq \frac{d}{y^l}.$$

PROOF. We prove the result for $l = 2$. The proof for the $l = 1$ case is similar and is left to the reader. From Proposition A.5, we know there exists a constant $c \in (0,\infty)$ such that $c^2 w_{k-1}(n+2) \geq w_{k-1}(n)$ for all $n \in \mathbb{Z}^+$ and all $k \in \{1,2,3\}$. The assumptions concerning $h$ imply the existence of an $M \in (0,\infty)$ such that $|(n+2)(n+1)h(n)| < M$ for all $n \in \mathbb{Z}^+$. Thus,

$$|E[h(N)w_{k-1}(N)|y]|$$
$$\leq E[|h(N)|w_{k-1}(N)|y]$$
$$= \sum_{n=0}^{\infty} |h(n)|w_{k-1}(n)\frac{e^{-y/2}(y/2)^n}{n!}$$
$$= \frac{4}{y^2} \sum_{n=0}^{\infty} (n+2)(n+1)|h(n)|w_{k-1}(n)\frac{e^{-y/2}(y/2)^{n+2}}{(n+2)!}$$
$$\leq \frac{4c^2}{y^2} \sum_{n=0}^{\infty} (n+2)(n+1)|h(n)|w_{k-1}(n+2)\frac{e^{-y/2}(y/2)^{n+2}}{(n+2)!}$$
$$\leq \frac{4c^2 M}{y^2} \sum_{n=0}^{\infty} w_{k-1}(n+2)\frac{e^{-y/2}(y/2)^{n+2}}{(n+2)!}$$
$$\leq \frac{4c^2 M}{y^2} \sum_{m=0}^{\infty} w_{k-1}(m)\frac{e^{-y/2}(y/2)^m}{m!}$$

$$= \frac{4c^2 M}{y^2} E[w_{k-1}(N)|y]. \qquad \square$$

PROPOSITION A.7. *Suppose $N|y \sim Poisson(y/2)$ where $y > 0$. For each $k \in \{1, 2, 3\}$, there exists a bounded function $\psi_k$ and a constant $d_k \in (0, \infty)$ such that $|\psi_k(y)| \leq d_k/y$ for $y > 0$ and*

$$\frac{E[w_k(N)|y]}{E[w_{k-1}(N)|y]} = \frac{p}{2} + 2(k - b - 1) + \frac{y}{2} + \psi_k(y).$$

PROOF. Suppose $a > 0$ and fix $k \in \{1, 2, 3\}$. The $a = 0$ case is simpler and is left to the reader. Rearranging (41) and taking expectations yields

$$E\left[\frac{w_{k-1}(N)}{N + p/2}\bigg|y\right]$$

$$= E\left[\frac{w_{k-1}(N+1)}{N + p/2 + k - b - 1 + \phi(N + k - 1)}\bigg|y\right]$$

$$= \frac{2}{y} \sum_{n=0}^{\infty} w_{k-1}(n+1)$$

(46)
$$\times \left[\frac{n+1}{n + p/2 + k - b - 1 + \phi(n + k - 1)}\right] \frac{e^{-y/2}(y/2)^{n+1}}{(n+1)!}$$

$$= \frac{2}{y} \sum_{m=1}^{\infty} w_{k-1}(m) \left[\frac{m}{m + p/2 + k - b - 2 + \phi(m + k - 2)}\right] \frac{e^{-y/2}(y/2)^m}{m!}$$

$$= \frac{2}{y} \sum_{m=0}^{\infty} w_{k-1}(m)[1 + O(m^{-1})] \frac{e^{-y/2}(y/2)^m}{m!}$$

$$= \frac{2}{y} E[w_{k-1}(N)|y] + \frac{2}{y} E[w_{k-1}(N) O(N^{-1})|y].$$

Therefore,

(47)
$$\frac{y(k - b - 1)}{2} \frac{E[w_{k-1}(N)/(N + p/2)|y]}{E[w_{k-1}(N)|y]}$$

$$= k - b - 1 + (k - b - 1) \frac{E[w_{k-1}(N) O(N^{-1})|y]}{E[w_{k-1}(N)|y]}.$$

This allows us to rewrite (43) as

$$\frac{E[w_k(N)|y]}{E[w_{k-1}(N)|y]} = \frac{p}{2} + 2(k - b - 1) + \frac{y}{2} + \psi_{k1}(y) + \psi_{k2}(y) + \psi_{k3}(y),$$

where
$$\psi_{k1}(y) := (k-b-1)\frac{E[w_{k-1}(N)O(N^{-1})|y]}{E[w_{k-1}(N)|y]},$$
$$\psi_{k2}(y) := \frac{y}{2}\frac{E[\phi(N+k-1)w_{k-1}(N)/(N+p/2)|y]}{E[w_{k-1}(N)|y]},$$
$$\psi_{k3}(y) := \frac{E[\phi(N+k-1)w_{k-1}(N)|y]}{E[w_{k-1}(N)|y]}.$$

Let $\psi_k(y) = \sum_{j=1}^{3}\psi_{kj}(y)$. Consider $\psi_{k1}(y)$. Proposition A.6 implies that there exists a $d_{k1}^* \in (0,\infty)$ such that
$$\left|\frac{E[w_{k-1}(N)O(N^{-1})|y]}{E[w_{k-1}(N)|y]}\right| \le \frac{d_{k1}^*}{y},$$
which implies that
$$|\psi_{k1}(y)| \le \frac{|k-b-1|d_{k1}^*}{y} = \frac{d_{k1}}{y}.$$

This shows that $\psi_{k1}$ is bounded for large $y$. Note that (47) yields
$$|\psi_{k1}(y)| = \left|b+1-k+\frac{y(k-b-1)}{2}\frac{E[w_{k-1}(N)/(N+p/2)|y]}{E[w_{k-1}(N)|y]}\right|$$
$$\le |b+1-k|+y\left|\frac{k-b-1}{p}\right|,$$
which shows that $\psi_{k1}$ is bounded for small $y$ as well.

Now consider $\psi_{k2}$. Since $\phi(n+k-1)/(n+p/2)$ is $O(n^{-2})$, Proposition A.6 implies that there exists a $d_{k2}^* \in (0,\infty)$ such that, for all $y > 0$,
$$\left|\frac{E[\phi(N+k-1)w_{k-1}(N)/(N+p/2)|y]}{E[w_{k-1}(N)|y]}\right| \le \frac{d_{k2}^*}{y^2},$$
and hence,
$$|\psi_{k2}(y)| \le \frac{y}{2}\frac{d_{k2}^*}{y^2} = \frac{d_{k2}}{y}.$$

As above, this shows that $\psi_{k2}$ is bounded for large $y$. The fact that $\psi_{k2}$ is bounded for small $y$ follows from the fact that $\phi$ is bounded.

Finally, consider $\psi_{k3}$. Since $\phi(n) = O(n^{-1})$, Proposition A.6 implies that there exists a $d_{k3}^* \in (0,\infty)$ such that, for all $y > 0$, $|\psi_{k3}(y)| < d_{k3}/y$. Again, the boundedness of $\psi_{k3}$ follows from that of $\phi$. Putting all of this together, we find that
$$|\psi_k(y)| = \left|\sum_{j=1}^{3}\psi_{kj}(y)\right| \le \sum_{j=1}^{3}|\psi_{kj}(y)| \le 3\frac{\max\{d_{k1},d_{k2},d_{k3}\}}{y} = \frac{d_k}{y}.$$

Moreover, since each $|\psi_{kj}(y)|$ is bounded for $y > 0$, so is $|\psi_k(y)|$. □

We can now assert that
$$\frac{E[w_1(N)|y]}{E[w_0(N)|y]} = \frac{p}{2} - 2b + \frac{y}{2} + \psi_1(y),$$
$$\frac{E[w_2(N)|y]}{E[w_1(N)|y]} = \frac{p}{2} + 2(1-b) + \frac{y}{2} + \psi_2(y),$$
$$\frac{E[w_3(N)|y]}{E[w_2(N)|y]} = \frac{p}{2} + 2(2-b) + \frac{y}{2} + \psi_3(y),$$
where each $\psi_i(y)$ is bounded and $|\psi_i(y)| \leq d/y$ for all $y > 0$ and a constant $d \in (0, \infty)$. We need one more technical result.

PROPOSITION A.8. *Suppose that $Y|\eta \sim \chi_p^2(\eta)$ and that $p \geq 3$. If there exists a $0 < d < \infty$ such that $|\psi(y)| \leq d/y$, then $|E[\psi(Y)|\eta]| \leq 2d/\eta$ so that $E[\psi(Y)|\eta] = O(\eta^{-1})$ as $\eta \to \infty$.*

PROOF.
$$|E[\psi(Y)|\eta]| \leq \int_0^\infty |\psi(y)| \sum_{n=0}^\infty \frac{e^{-\eta/2}(\eta/2)^n}{n!} \frac{(y/2)^{n+p/2-1}e^{-y/2}}{2\Gamma(n+p/2)} dy$$
$$\leq d \sum_{n=0}^\infty \frac{e^{-\eta/2}(\eta/2)^n}{n!} \int_0^\infty \frac{y^{n+p/2-2}e^{-y/2}}{2^{n+p/2}\Gamma(n+p/2)} dy$$
$$= d \sum_{n=0}^\infty \frac{e^{-\eta/2}(\eta/2)^n}{n!} \frac{1}{(2n+p-2)}$$
$$= \frac{2d}{\eta} \sum_{n=0}^\infty \frac{(n+1)}{(2n+p-2)} \frac{e^{-\eta/2}(\eta/2)^{n+1}}{(n+1)!}$$
$$\leq \frac{2d}{\eta} \sum_{n=0}^\infty \frac{e^{-\eta/2}(\eta/2)^{n+1}}{(n+1)!} \leq \frac{2d}{\eta}.$$
□

Assume that $Y|\eta \sim \chi_p^2(\eta)$. Proposition A.8 implies that $E[|\psi_i(Y)||\eta] = O(\eta^{-1})$ as $\eta \to \infty$. Furthermore, since $|y\psi_i(y)| \leq d$, $E[|Y\psi_i(Y)||\eta] = O(1)$ as $\eta \to \infty$. Also, since the $\psi_i$ are bounded, so are products of the $\psi_i$. Thus, $E[|\psi_i(Y)\psi_j(Y)||\eta]$ and $E[|\psi_i(Y)\psi_j(Y)\psi_k(Y)||\eta]$ are both $O(1)$ as $\eta \to \infty$ for any $i, j, k \in \{1, 2, 3\}$. Finally, since $|y^2 \psi_i(y)| \leq dy$, it follows that $E[Y^2|\psi_i(Y)||\eta] \leq dE[Y|\eta] = d(\eta + p) = O(\eta)$ as $\eta \to \infty$.

Putting the above work together (and using the moments of the noncentral $\chi^2$ given at the beginning of this section), we calculate that
$$E\left[\frac{E(w_1(N)|Y)}{E(w_0(N)|Y)}\bigg|\eta\right] = \frac{\eta}{2} + p - 2b + O(\eta^{-1})$$

(in all of these equations, it is understood that the limits are taken as $\eta \to \infty$)

$$E\left[\frac{E(w_2(N)|Y)}{E(w_0(N)|Y)}\Big|\eta\right] = E\left[\frac{E(w_2(N)|Y)}{E(w_1(N)|Y)}\frac{E(w_1(N)|Y)}{E(w_0(N)|Y)}\Big|\eta\right]$$

$$= \frac{\eta^2}{4} + (p+2-2b)\eta + O(1),$$

$$E\left[\frac{E(w_3(N)|Y)}{E(w_0(N)|Y)}\Big|\eta\right] = E\left[\frac{E(w_3(N)|Y)}{E(w_2(N)|Y)}\frac{E(w_2(N)|Y)}{E(w_1(N)|Y)}\frac{E(w_1(N)|Y)}{E(w_0(N)|Y)}\Big|\eta\right]$$

$$= \frac{1}{8}(\eta^3 + (6p - 12b + 24)\eta^2) + O(\eta).$$

We are now in position to calculate $\mu_1(\eta)$, $\mu_2(\eta)$ and $\mu_3(\eta)$ in (38):

$$\mu_1(\eta) = 2E\left[\frac{E(w_1(N)|Y)}{E(w_0(N)|Y)}\Big|\eta\right] - \eta$$

$$= \eta + 2p - 4b + O(\eta^{-1}) - \eta$$

$$= 2p - 4b + O(\eta^{-1}),$$

$$\mu_2(\eta) = 4E\left[\frac{E(w_2(N)|Y)}{E(w_0(N)|Y)}\Big|\eta\right] - 4\eta E\left[\frac{E(w_1(N)|Y)}{E(w_0(N)|Y)}\Big|\eta\right] + \eta^2$$

$$= 4\left(\frac{\eta^2}{4} + (p+2-2b)\eta\right) - 4\eta\left(\frac{\eta}{2} + p - 2b\right) + \eta^2 + O(1)$$

$$= 8\eta + O(1),$$

$$\mu_3(\eta) = 8E\left[\frac{E(w_3(N)|Y)}{E(w_0(N)|Y)}\Big|\eta\right] - 12\eta E\left[\frac{E(w_2(N)|Y)}{E(w_0(N)|Y)}\Big|\eta\right]$$

$$+ 6\eta^2 E\left[\frac{E(w_1(N)|Y)}{E(w_0(N)|Y)}\Big|\eta\right] - \eta^3$$

$$= (\eta^3 + (6p - 12b + 24)\eta^2) - 12\eta\left(\frac{\eta^2}{4} + (p+2-2b)\eta\right)$$

$$+ 6\eta^2\left(\frac{\eta}{2} + p - 2b\right) - \eta^3 + O(\eta) = O(\eta).$$

This completes the proof of Proposition A.1. □

**Acknowledgments.** The authors are grateful to Brian Shea and an anonymous referee for insightful comments that led to substantial improvements in the paper.

M. L. EATON
G. L. JONES
SCHOOL OF STATISTICS
UNIVERSITY OF MINNESOTA
MINNEAPOLIS, MINNESOTA 55455
USA
E-MAIL: eaton@stat.umn.edu
　　　　galin@stat.umn.edu

J. P. HOBERT
DEPARTMENT OF STATISTICS
UNIVERSITY OF FLORIDA
GAINESVILLE, FLORIDA 32611
USA
E-MAIL: jhobert@stat.ufl.edu

W.-L. Lai
Department of International Business
Providence University
Taichung
Taiwan
E-mail: wllai@pu.edu.tw